\documentclass{article}

\usepackage{amssymb} 
\usepackage{amsmath} 
\usepackage{latexsym} 
\usepackage{theorem} 
\usepackage{color}
\usepackage{epsfig}

\theorembodyfont{\itshape} 

\newtheorem{theorem}{Theorem}[section]
\newtheorem{lemma}[theorem]{Lemma}

\newtheorem{corollary}[theorem]{Corollary}

\newtheorem{definition}[theorem]{Definition}

\newtheorem{remark}[theorem]{Remark}
\newenvironment{proof}{{\par\addvspace{0.1cm}\noindent \bf Proof. }}{\hfill$\Box$\par\medskip} 
 
\def\n{m}
\def\a{\alpha}

\def\an{{\alpha-\n}}
\def\d{\delta}

\def\Om{\varOmega}
\def\Ga{\varGamma}
\def\RR{\mathbb{R}}

\def\diam{{d}}
\def\rhosub#1{\mbox{\large $\rho $}_{\mbox{\small $\!{}_{#1}$}}}

\def\vect#1{\mbox{\boldmath $#1$}} 
\def\interior#1{\stackrel{\circ}{#1}}

\title{\bf Uniqueness of radial centers of parallel bodies}
\author{Jun O'Hara}

\numberwithin{equation}{section}

\begin{document}

\maketitle

\begin{abstract} 
We show the uniqueness of the radial centers of any order $\a$ of a parallel body of a convex body $\Om$ in $\RR^\n$ at distance $\d$ 
if $\d$ is greater than the diameter of $\Om$ multiplied by a constant which depends only on the dimension $\n$. 
\end{abstract}



\section{Introduction}
Let $\Om$ be a {\em body} in $\RR^\n$ $(\n\ge2)$, i.e. a compact set which is a closure of its interior, with a piecewise $C^1$ boundary. 
Consider a potential of the form 
\[
V_\Om^{(\a)}(x)=\int_\Om{|x-y|}^{\an}\,d\mu(y) \hspace{0.5cm}(\a>0),\]
where $\mu$ is the standard Lesbegue measure of $\RR^\n$. 
It is a singular integral when $\a<\n$ and $x\in\Om$. 
When $0<\a<m$ it is the {\em Riesz potential} of (the characteristic function $\chi_\Om$ of) $\Om$. 

In particular, when $\Om$ is convex and $x\in\interior\Om$,  $V_\Om^{(\a)}(x)$ can be expressed as 
$$
V_\Om^{(\a)}(x)=\frac1\a\int_{S^{\n-1}}\left(\rhosub{\Om_{-x}}(v)\right)^\a\,d\sigma(v) 
$$
where $\sigma$ is the standard Lebesgue measure of $S^{\n-1}$ and $\rhosub{\Om_{-x}}:S^{\n-1}\to\RR_{>0}$ is a {\em radial function} of $\Om_{-x}=\{y-x\,|\,y\in\Om\}$ given by $\mbox{\large $\rho  $}_{\Om_{-x}}(v)=\sup\{a\ge0\,|\,x+av\in \Om\}$. 
Thus $V_\Om^{(\a)}(x)$ coincides with the {\em dual mixed volume} 
as introduced by Lutwak (\cite{Lu75,Lu88}) up to multiplication by a constant. 

In \cite{O} we defined an {\em $r^{\an}$-center of $\Om$}. It is a point where the extreme value of $V_\Om^{(\a)}$ (minimum or maximum according to the value of $\a$) is attained when $\a\ne\n$. (The case when $\a=\n$ will be addressed later.) 
For example, the center of mass is an $r^2$-center. 
When $\Om$ is convex, an $r^{\an}$-center $(\a\ne\n)$ coincides with the {\em radial center of order $\a$}, which was introduced in \cite{M1} for $0<\a\le1$ and in \cite{H} in general $(\a\ne0)$. 
An $r^{\an}$-center of a body $\Om$ exists for any $\a$ and is unique 
if $\a\ge\n+1$ or if $\a\le1$ and $\Om$ is convex (\cite{O}). 
It was conjectured that a convex subset $\Om$ has a unique $r^{\an}$-center for any $\a$. 

In this paper we show the uniqueness of an $r^\an$-center for any $\a$ when $\Om$ is close to a ball in some sense. 
To be precise, we show that there is a positive function $\varphi(\n)$ such that for any convex body $\widetilde\Om$ with a piecewise $C^1$ boundary, if $\d\ge\varphi(\n)\cdot\textrm{diam}(\widetilde\Om)$ then a $\d$-{\em parallel body} of $\widetilde\Om$ has a unique $r^\an$-center for any $\a$. 
Here, a $\d$-parallel body of $\widetilde\Om$ is the closure of a $\d$-tubular neighbouhood of $\widetilde\Om$, and is denoted by $\widetilde\Om+\d B^\n$. 
%
%
The proof has two steps. 
First we show that a center can appear only in $\widetilde\Om$ 
by the so-called {\sl moving plane method} in analysis (\cite{GNN}). 
Then we show that $V_{\widetilde\Om+\d B^\n}^{(\a)}$ is convex (or concave according to the value of $\a$) on $\widetilde\Om$ using the boundary integral expression of the second derivatives of $V_\Om^{(\a)}$. 

\medskip
%
Throughout the paper, $\interior X$, $\overline{X}$, $X^c$, and $\textrm{conv}(X)$ denote the interior, the closure, the complement, and the convex hull of $X$ respectively. 
We denote the standard Lesbegue measure of $\RR^\n$ by $\mu$, and that of $\partial\Om$ and other $(\n-1)$-dimensional spaces like $S^{\n-1}$ by $\sigma$.

\section{Preliminaries from \cite{O}}
%
In this section we introduce some of the results of \cite{O} which are necessary for our study. 

First remark that if we define 
$$X-Y=(X\setminus(X\cap Y))\cup-(Y\setminus(X\cap Y)) \hspace{0.5cm}(X,Y\subset \RR^\n),$$
where the second term is equipped with the reverse orientation, 
then 
\[V^{(\a)}_{\Om_1-\Om_2}(x)=V^{(\a)}_{\Om_1}(x)-V^{(\a)}_{\Om_2}(x) \hspace{0.5cm}(x\in\interior{\Om}_1\cap\interior{\Om}_2)\] 
for any $\a$. 

%
%
\subsection{Boundary integral expression of the derivatives}
The first derivatives of $V_{\Om}^{(\a)}$ can be expressed by the boundary integral as 
\begin{equation}\label{formula_derivative_boundary}
\frac{\partial V_{\Om}^{(\a)}}{\partial x_j}(x)=-\int_{\partial\Om}{|x-y|}^{\an}\,e_j\cdot n\,d\sigma (y)
\end{equation}
for any $j$ $(1\le j\le \n)$ if $x=(x_1, \ldots, x_\n)\not\in\partial\Om$, where $n$ is a unit outer normal vector to $\partial\Om$ at $y$, $e_j$ is the $j$-th unit vector of $\RR^\n$, and $\sigma$ denotes the standard Lebesgue measure of $\partial\Om$. 
%
%
This is because 
$$
\frac{\partial \, r^{\an}}{\partial x_j}=-\frac{\partial \, r^{\an}}{\partial y_j}=-\textrm{div}_y\left(r^\an \,e_j\right). 
$$
It follows that the second derivatives satisfy 
\begin{equation}\label{f_second_partial_derivative_boundary}
\displaystyle \frac{\partial^2 V_{\Om}^{(\a)}}{\partial x_j{}^2}(x)
=\displaystyle -(\a-\n)\int_{\partial\Om}{|x-y|}^{\an-2}(x_j-y_j)\,e_j\cdot n\,d\sigma (y)  
\end{equation}
for any $\a$ if $x\not\in\partial\Om$ (or for any $x$ if $\a>2$). 
Furthermore, if $x\in\Om^c$ then for any $\a$ 
\begin{eqnarray}
\displaystyle \frac{\partial^2 V_{\Om}^{(\a)}}{\partial x_j{}^2}(x)
%
&=&\displaystyle (\a-\n)\int_{\Om}{|x-y|}^{\an-4}\left((\an-2)(x_j-y_j)^2+{|x-y|}^2\right)d\mu (y) \label{f_second_partial_derivative_Omega_bis}\\
&=&\displaystyle (\a-\n)\int_{\Om}{|x-y|}^{\an-4}\left((\an-1)(x_j-y_j)^2+\sum_{i\ne j}(x_i-y_i)^2\right)\!d\mu (y)\,.\hspace{0.3cm}{\phantom{a}}
\label{f_second_partial_derivative_Omega}
\end{eqnarray}
%

%
\subsection{Definition of the $\vect{r^\an}$-centers}
When $\a\ne\n$ we call a point {\em ${r^{\an}}$-center} of $\Om$ if it gives the minimum value of $V_{\Om}^{(\a)}$ when $\a>\n$ and the maximum value of $V_{\Om}^{(\a)}$ when $0<\a<\n$. 
When $\a=\n$ it is meaningless to use $V_{\Om}^{(\n)}$ as it is constantly equal to $\textrm{Vol}\,(\Om)$. 
We call a point an {\em $r^0$-center} if it gives the maximum value of the log potential \[
V_\Om^{\log}(x)=\int_\Om\log\frac1{|x-y|}\,d\mu(y)=-\int_\Om\log{|x-y|}\,d\mu(y).
\]
%

As we noticed in the introduction, the center of mass is an $r^2$-center, and if $\Om$ is convex and $\a\ne\n$, an $r^{\an}$-center coincides with the {\em radial center of order $\a$}, which was introduced in \cite{M1} for $0<\a\le1$ and in \cite{H} for $\a\ne0$. 

We remark that the statements of $r^\an$-centers in the case when $\a=\n$ in this paper can be obtained exactly in the same way as in the case when $0<\a<\n$. 
This is because we only use the estimate on the second derivative in our study, and that of the log potential  
\[\begin{array}{rcl}
\displaystyle \frac{\partial^2 V_\Om^{\log}}{\partial {x_j}^2}(x)&=&\displaystyle \int_{\partial\Om}{|x-y|}^{-2}(x_j-y_j)e_j\cdot n\,d\sigma (y)
\end{array}
\]
can be considered as the limit of $1/(\n-\a)$ times the second derivative of $V_{\Om}^{(\a)}$ as $\a$ approaches $\n$ (see \eqref{f_second_partial_derivative_boundary}). 
%
\subsection{Minimal unfolded regions}\label{subsection_minimal_unfolded_region}
%
Let $v$ be a unit vector in $S^{\n-1}$ and $b$ be a real number. Put 
$$
H_{v,b}=\{x\in\RR^\n\,|\,x\cdot v=b\}, \> H_{v,b}^+=\{x\in\RR^\n\,|\,x\cdot v>b\}, \> H_{v,b}^-=\{x\in\RR^\n\,|\,x\cdot v<b\}.
$$
Let $\textrm{Refl}_{v,b}$ be a reflection of $\RR^\n$ in $H_{v,b}$. 
Let $\Om$ be a compact set in $\RR^\n$. 
Put $\displaystyle M_v=M_v(\Om)=\max_{x\in\Omega}x\cdot v$ 
and 
$$
u_v=u_v(\Om)=\inf\Big\{a\,\big|\, a\le M_v, \,\textrm{Refl}_{v,b}\big(\Om\cap H_{v,b}^+\big)\subset\Om
\>\>(a\le{}^{\mbox{$\forall$}}b\le M_v)\Big\}.
$$
%
%
The {\em minimal unfolded region} of $\Om$ is given by 
$$
Uf(\Om)
=\bigcap_{v\in S^{\n-1}} \overline{H_{v,u_v}^-}.
$$
It is a non-empty compact convex set and is contained in the convex hull of $\Om$. 
It is not necessarily contained in $\Om$. 
\begin{figure}[htbp]
\begin{center}
\begin{minipage}{.4\linewidth}
\begin{center}
\includegraphics[width=.6\linewidth]{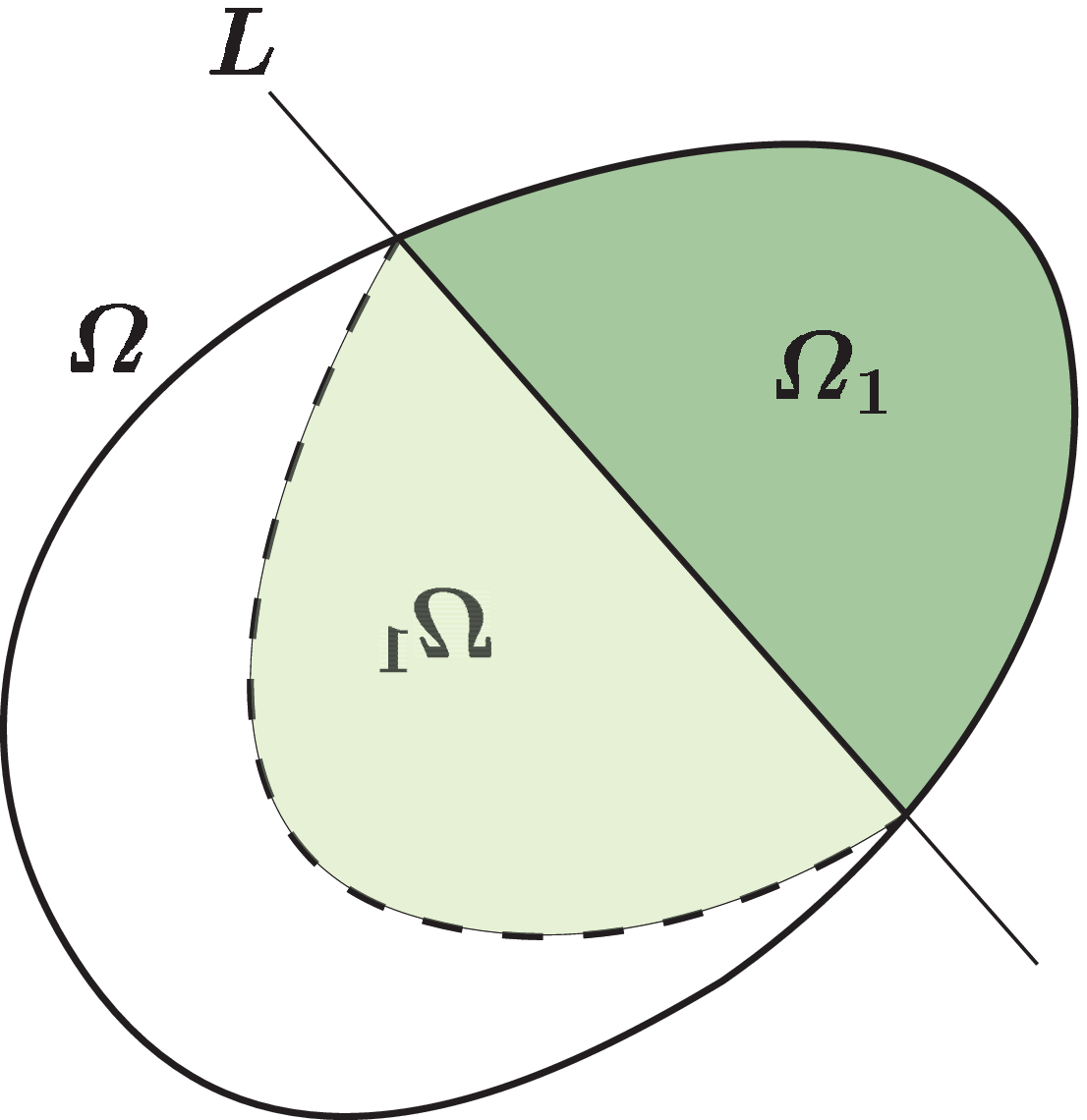}
\caption{Folding a convex set like origami}
\label{origami_convex}
\end{center}
\end{minipage}
\hskip 0.4cm
\begin{minipage}{.55\linewidth}
\begin{center}
\includegraphics[width=.66\linewidth]{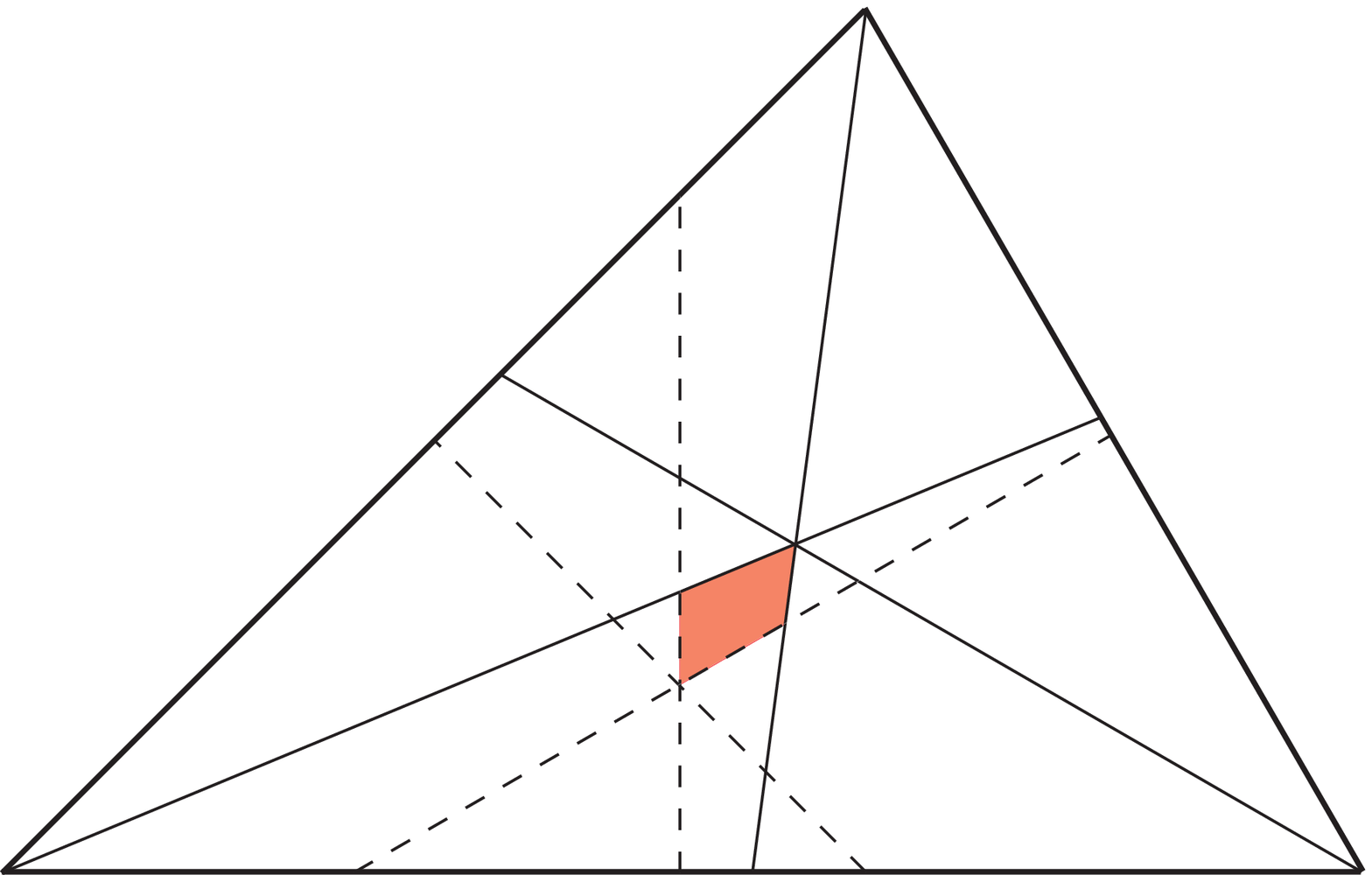}
\caption{A minimal unfolded region of a non-obtuse triangle. Bold lines are angle bisectors, dotted lines are perpendicular bisectors. }
\label{folded_core_triangle}
\end{center}
\end{minipage}
\end{center}
\end{figure}

\subsection{Existence and uniqueness of $\vect{r^\an}$-centers}
%
\begin{theorem}{\rm (\cite{O})}\label{thm_O} 
Let $\Om$ be a body in $\RR^\n$ with a piecewise $C^1$ boundary $\partial\Om$. 
\begin{enumerate}
\item There exists an ${r^{\an}}$-center of $\Om$ for any $\a$. 
\item An ${r^{\an}}$-center is contained in the minimal unfolded region of $\Om$ for any $\a$. 
\item An ${r^{\an}}$-center of $\Om$ is unique if $\a\ge \n+1$. 
\item An ${r^{\an}}$-center of $\Om$ is unique if $\a\le1$ and $\Om$ is convex. 
\end{enumerate}
\end{theorem}
The second statement is essentially based on the so-called {\sl moving plane method} in analysis (\cite{GNN}) as the integrands appearing in the formulae of $V_\Om^{(\a)}$ and its derivatives are symmetric. 
The uniqueness of an ${r^{\an}}$-center follows from $\displaystyle \frac{\partial^2 V_{\Om}^{(\a)}}{\partial {x_j}^2}>0$ when $\a\ge \n+1$ and the strong concavity $V_{\Om}^{(\a)}$ on $\interior\Om$ when $\a\le1$ and $\Om$ is convex.

\section{Uniqueness of the centers of $\vect{\Om+\d B^\n}$}
%
We conjectured that if $\Om$ is convex then it has only one ${r^{\an}}$-center for any $\a$, although it was proved that $V_{\Om}^{(\a)}$ is not necessarily convex nor concave. 
In this section we show that the conjeture holds for a $\delta$-parallel bodies $\widetilde\Om+\d B^\n\{x+u\,|\,x\in\widetilde\Om, \,u\in B^\n\}$ provided that $\delta$ is large enough compared with the diameter of $\widetilde\Om$.  
To be precise, we prove the following theorem: 
\begin{theorem}\label{main_thm1}
For any natural number $\n\ge2$ there is a positive constant $\varphi(\n)$ such that for any compact convex set $\widetilde\Om$ in $\RR^\n$ with piecewise $C^1$ boundary, if $\d\ge\varphi(\n)\cdot\textrm{\rm diam}(\widetilde\Om)$ then $\widetilde\Om+\d B^\n$ has a unique $r^\an$-center for any $\a$. 
\end{theorem}
%

By Theorem \ref{thm_O} 
it is enough to show 
\begin{equation}\label{f_second_der}
\displaystyle \frac{\partial^2 V_{\widetilde\Om+\d B^\n}^{(\a)}}{\partial {x_j}^2}<0
 \>\>\>(1<\a<\n), 
\>\frac{\partial^2 V_{\widetilde\Om+\d B^\n}^{\log}}{\partial {x_j}^2}<0, 
\>\frac{\partial^2 V_{\widetilde\Om+\d B^\n}^{(\a)}}{\partial {x_j}^2}>0\>\>\>(\n<\a<\n+1)
\end{equation}
%
for any $j$ on the minimal unfolded region of $\widetilde\Om+\d B^\n$. 
\begin{lemma}\label{lemma_uf_Omega+deltaB}
Let $\widetilde\Om$ be any compact subset of $\RR^\n$. 
The minimal unfolded region of $\widetilde\Om+\d B^\n$ is contained in the convex hull of $\widetilde\Om$ for any $\d>0$. 
\end{lemma}
\begin{proof}
Let us use the notation in Subsection \ref{subsection_minimal_unfolded_region}. 
Let $v\in S^{\n-1}$ be any vector and $b$ any real number that satisfies
$$
M_v(\widetilde\Om)<b\le M_v(\widetilde\Om+\d B^\n)=M_v(\widetilde\Om)+\d.
$$
Then for any point $Q$ in $\widetilde\Om$ we have $\textrm{Refl}_{v,b}\big(B_\d(Q)\cap H_{v,b}^+\big)\subset B_\d(Q)\cap H_{v,b}^-$ as the center $Q$ is in $H_{v,b}^-$. 
As $(\widetilde\Om+\d B^\n)\cap H_{v,b}^+=\bigcup_{Q\in\widetilde\Om}\big(B_\d(Q)\cap H_{v,b}^+\big)$ we have 
$$\begin{array}{rcl}
\textrm{Refl}_{v,b}\big((\widetilde\Om+\d B^\n)\cap H_{v,b}^+\big)
&=&\displaystyle \bigcup_{Q\in\widetilde\Om}\textrm{Refl}_{v,b}\big(B_\d(Q)\cap H_{v,b}^+\big)\\[6mm]
&\subset&\displaystyle \bigcup_{Q\in\widetilde\Om} \big(B_\d(Q)\cap H_{v,b}^-\big)\\[6mm]
&\subset&\displaystyle (\widetilde\Om+\d B^\n)\cap H_{v,b}^-. 
\end{array}$$
Consequently we have $u_v(\widetilde\Om+\d B^\n)\le M_v(\widetilde\Om).$ 
It follows that 
$$
Uf(\widetilde\Om+\d B^\n)=\bigcap_{v\in S^{\n-1}} \overline{H_{v,u_v(\widetilde\Om+\d B^\n)}^-}
\subset \bigcap_{v\in S^{\n-1}} \overline{H_{v,M_v(\widetilde\Om)}^-}=\textrm{conv}(\widetilde\Om).
$$
\end{proof}

Next we proceed to the proof of \eqref{f_second_der} on $\widetilde\Om$. 

\begin{definition} \rm 
Let $\n$ be a natural number with $\n\ge2$ and let $a>0$. 
Let $L_a$ denote an oriented line segment in $\RR^2$ which starts from $(a,0)$ and ends at $(0,1)$. 
For real numbers $\a$ and $\xi$ with $0\le\xi< a$ define 
\begin{equation}\label{def_F}
F(\n,\a,a,\xi)=\int_{L_a}{|x-y|}^{\an-2}(\xi-y_1)\,{y_2}^{\n-2}\,dy_2\,,
\end{equation}
where $x=(\xi,0)$. 
\end{definition}
\begin{lemma}\label{key_lemma_m=2} 
Suppose $\n=2$ and $1<\a<3$. For any $a>0$, if $0\le \xi\le\frac{a}2$ then $F(2,\a,a,\xi)<0$. 
%
\end{lemma}
\begin{proof}
Divide $L_a$ into three parts;
$$
L_1=L_a\cap\{2\xi\le x_1\le a\},\,
L_2=L_a\cap\{\xi\le x_1\le 2\xi\},\,
L_3=L_a\cap\{0\le x_1\le \xi\}.
$$
If we put 
\begin{figure}[htbp]
\begin{center}
\includegraphics[width=.36\linewidth]{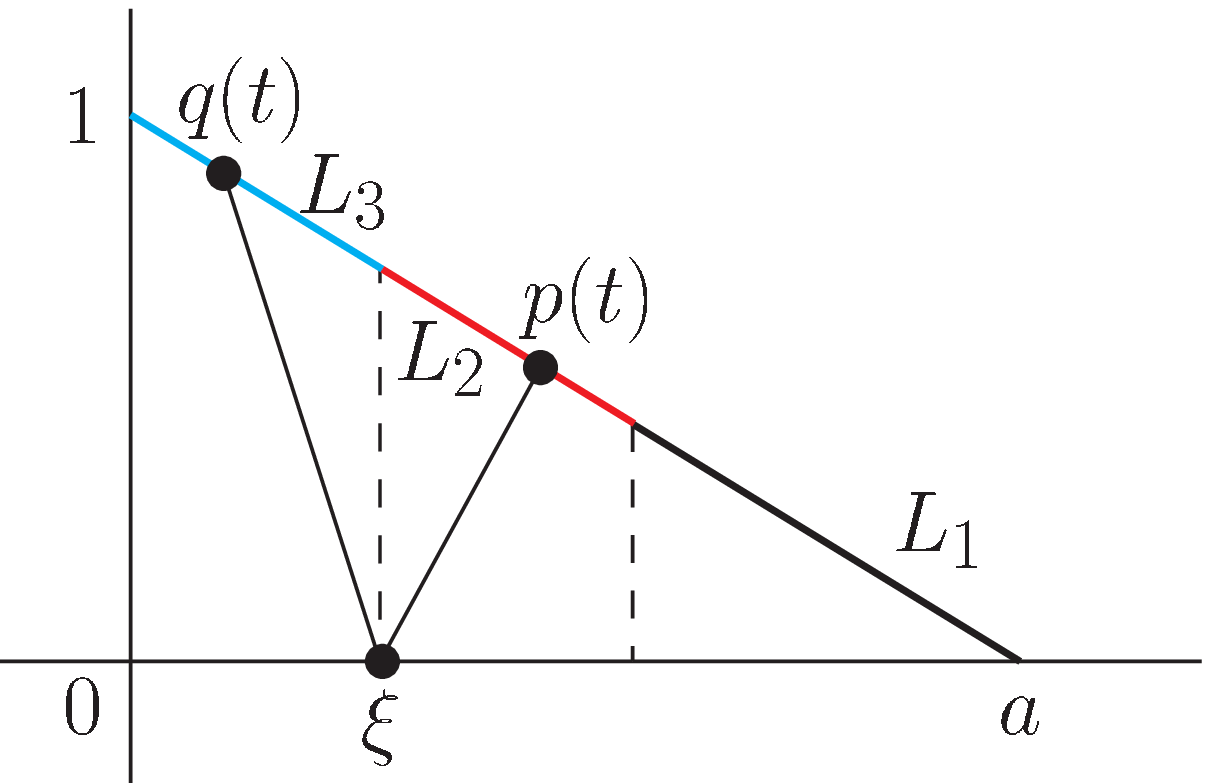}
\caption{}
\label{L1L2L3}
\end{center}
\end{figure}
$$
p(t)=\left(\xi+t,\left(1-\frac\xi{a}\right)-\frac{t}a\right)\in L_2, \,
q(t)=\left(\xi-t,\left(1-\frac\xi{a}\right)+\frac{t}a\right)\in  L_3 \hspace{0.5cm}(0\le t\le\xi) 
$$
then $|x-p(t)|<|x-q(t)|$ (Figure \ref{L1L2L3}). Hence, as $y_1=a(1-y_2)$ on $L_a$, 
$$
\int_{L_2\cup L_3}{|x-y|}^{\a-4}(\xi-y_1)\,dy_2
=\int_0^\xi\left(-{|x-p(t)|}^{\a-4}+{|x-q(t)|}^{\a-4}\right)\cdot t\cdot\frac{dt}{a}
<0.$$
Since
$\int_{L_1}{|x-y|}^{\a-4}(\xi-y_1)\,dy_2<0$, it completes the proof. 
\end{proof}
%
%
\begin{lemma}\label{lemma_c_a_m>2}
Suppose $\n\ge3$ and $1<\a<\n+1$. 
There is a map $\psi_\a\colon \RR_+\to\RR_+$ such that for any $c>0$ if $a\ge \psi_\a(c)\,c$ then 
\begin{equation}\label{f_-a^c}
\int_{-a}^c\frac{(t-\frac1a){(t+a)}^{\n-2}}{{\left(\sqrt{t^2+1}\,\right)}^{\n+2-\a}}\,dt<0.
\end{equation}
\end{lemma}
\begin{proof}
Let $c>0$. 
Assume $a> 2c$. Note that the integrand of \eqref{f_-a^c} is positive (or negative) if $t>\frac1a$ (or respectively $t<\frac1a$). 
Therefore we have only to consider the case when $\frac1a<c$. 
Observe that 
\begin{equation}\label{f_ineq_int_-a^c}
\int_{-a}^c\frac{(t-\frac1a){(t+a)}^{\n-2}}{{\left(\sqrt{t^2+1}\,\right)}^{\n+2-\a}}\,dt\le\displaystyle 
-\int_{[-2c,-c]\cup[-c,-\frac1a]}\frac{|t-\frac1a|{(t+a)}^{\n-2}}{{\left(\sqrt{t^2+1}\,\right)}^{\n+2-\a}}\,dt
+\int_{\frac1a}^{c}\frac{(t-\frac1a){(t+a)}^{\n-2}}{{\left(\sqrt{t^2+1}\,\right)}^{\n+2-\a}}\,dt,
%
\end{equation}
where the right hand side can be estimated by
$$
\begin{array}{rcl}
\displaystyle \int_{-c}^{-\frac1a}\frac{|t-\frac1a|{(t+a)}^{\n-2}}{{\left(\sqrt{t^2+1}\,\right)}^{\n+2-\a}}\,dt
&\ge&\displaystyle \left(\frac{a-c}{a+c}\right)^{\n-2}\int_{\frac1a}^{c}\frac{(t-\frac1a){(t+a)}^{\n-2}}{{\left(\sqrt{t^2+1}\,\right)}^{\n+2-\a}}\,dt\,, \\[5mm]
\displaystyle \int_{-2c}^{-c}\frac{|t-\frac1a|{(t+a)}^{\n-2}}{{\left(\sqrt{t^2+1}\,\right)}^{\n+2-\a}}\,dt
&\ge&\displaystyle \left(\frac{a-2c}{a+c}\right)^{\n-2}\cdot\frac1{{\left(\sqrt{4c^2+1}\,\right)}^{\n+2-\a}}\int_{\frac1a}^{c}\frac{(t-\frac1a){(t+a)}^{\n-2}}{{\left(\sqrt{t^2+1}\,\right)}^{\n+2-\a}}\,dt\,. \\[4mm]
\end{array}
$$
As 
\begin{equation}\label{ineq1}
\left(\frac{a-c}{a+c}\right)^{\n-2}+\left(\frac{a-2c}{a+c}\right)^{\n-2}\frac1{{\left(\sqrt{4c^2+1}\,\right)}^{\n+2-\a}}
\ge \left(\frac{a-2c}{a+c}\right)^{\n-2}\left(1+{(4c^2+1)}^{-\frac{\n+2-\a}2}\right), 
\end{equation}
if we put 
\begin{equation}\label{def_psi}
\psi_\a(c)=\frac{2{\left(1+{(4c^2+1)}^{-\frac{\n+2-\a}2}\right)}^{\frac1{\n-2}}+1}{{\left(1+{(4c^2+1)}^{-\frac{\n+2-\a}2}\right)}^{\frac1{\n-2}}-1}
=2+\frac3{{\left(1+{(4c^2+1)}^{-\frac{\n+2-\a}2}\right)}^{\frac1{\n-2}}-1}
\end{equation}
then $a\ge\psi_\a(c)c$ is equivalent to 
$$
\frac{a-2c}{a+c}\ge{\left(1+{(4c^2+1)}^{-\frac{\n+2-\a}2}\right)}^{-\frac1{\n-2}}\,,
$$
which implies that the right hand side of \eqref{f_ineq_int_-a^c} is negative; thus the proof is completed. 
Remark that as $\psi_\a(c)>2$, if $a\ge \psi_\a(c)\,c$ then $a$ satisfies the assumption $a>2c$ which appeared at the beginning of the proof. 
%
%
\end{proof}
%
%
\begin{corollary}\label{cor_lemma_c_a_m>2}
Suppose $\n\ge3$ and $1<\a<\n+1$. 
For any $\xi_0>0$ there is $a_0>0$ such that if $0\le\xi\le\xi_0$ and $a\ge a_0$ then 
$F(\n,\a,a,\xi)<0$, where $F(\n,\a,a,\xi)$ is given by \eqref{def_F}. 
In fact, we can take 
$$
a_0=\psi_\a\left(\frac{2\xi_0^2+1}{\xi_0}\right)\frac{2\xi_0^2+1}{\xi_0},
$$
where $\psi_\a$ is given by \eqref{def_psi}. 
\end{corollary}

\begin{figure}[htbp]
\begin{center}
\includegraphics[width=.4\linewidth]{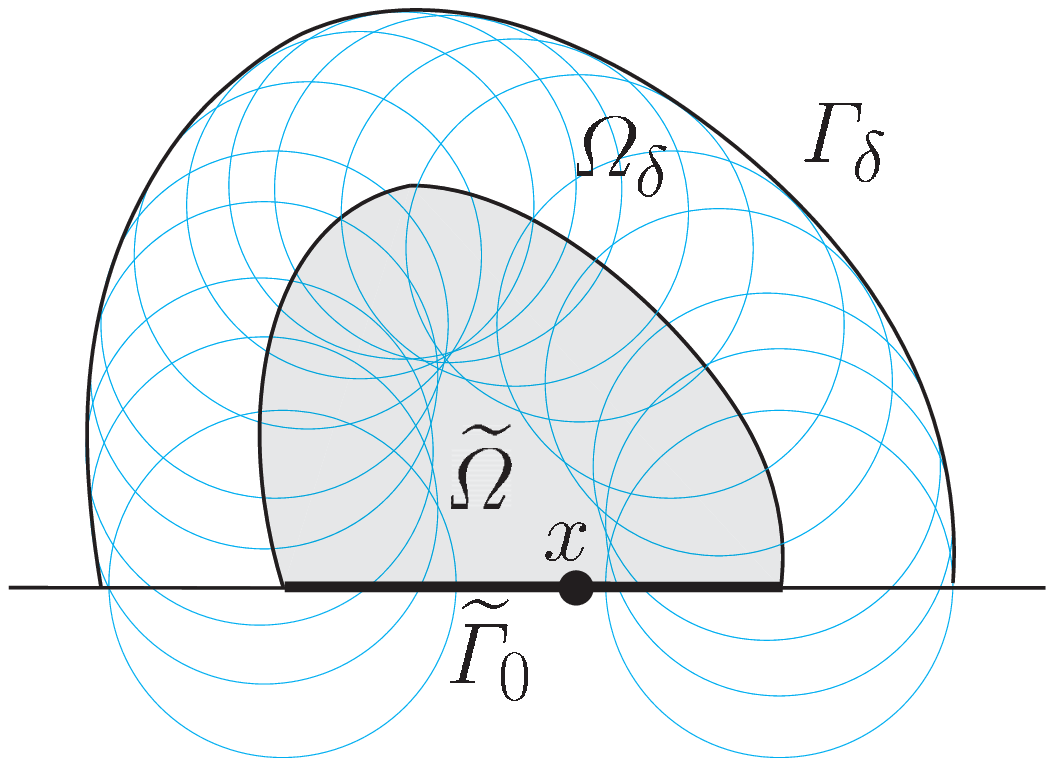}
\caption{$\d$-parallel body $\Om=\Om_\d$. $\Ga_\d$ is an envelope of circles with radius $\d$ whose centers lie on $\partial\widetilde\Om\cap\RR^2_+$}
\label{parallel_body}
\end{center}
\end{figure}

\begin{proof} 
As $y_1=a(1-y_2)$ on $L_a$,
\begin{eqnarray}\label{f_key_lemma_target}
F(\n,\a,a,\xi)=\int_0^1\left(\left(a(1-y_2)-\xi\right)^2+{y_2}^2\right)^{\frac\an2-1}(\xi-a(1-y_2))\,{y_2}^{\n-2}\,dy_2.
\end{eqnarray}
Put $\displaystyle y_2-\frac{a(a-\xi)}{1+a^2}=\frac{a-\xi}{1+a^2}\,t$. 
Then, as 
\[\begin{array}{l}
\displaystyle 
\left(a(1-y_2)-\xi\right)^2+{y_2}^2=\frac{(a-\xi)^2}{1+a^2}(t^2+1), \> 
\displaystyle \xi-a(1-y_2)=\frac{a-\xi}{1+a^2}(at-1), \> 
\displaystyle y_2=\frac{a-\xi}{1+a^2}(t+a),
\end{array}\]
%
$F(\n,\a,a,\xi)<0$ is equivalent to 
\[
\int_{-a}^{\frac{a\xi+1}{a-\xi}}
\frac{(t-\frac1a){(t+a)}^{\n-2}}{{\left(\sqrt{t^2+1}\,\right)}^{\n+2-\a}}\,dt<0.
\]

First remark that if $0<\xi<\xi_0$ then $F(\n,\a,a,\xi)<F(\n,\a,a,\xi_0)$ since $\frac{a\xi+1}{a-\xi}$ is an increasing function of $\xi$ $(0<\xi<a)$ with $\lim_{\xi\to0}\frac{a\xi+1}{a-\xi}=\frac1a$ and the integrand is positive when $t>\frac1a$. 

On the other hand, if we put $c(a)=\frac{a\xi_0+1}{a-\xi_0}$, it is a decreasing function of $a$ as $c(a)=\xi_0+\frac{1+\xi_0^2}{a-\xi_0}$. 
Put 
$$
c_0=c(2\xi_0)=\frac{2\xi_0^2+1}{\xi_0}, \>\>
a_0=\psi_\a(c_0)c_0=\psi_\a\left(\frac{2\xi_0^2+1}{\xi_0}\right)\frac{2\xi_0^2+1}{\xi_0},
$$
where $\psi_\a$ is given by \eqref{def_psi}. 
If $a\ge a_0$ then $c(a)<c(2\xi_0)=c_0$ as $a>2\xi_0$. Since $c(a)>\frac1a$ we have 
$$
\int_{-a}^{c(a)}
\frac{(t-\frac1a){(t+a)}^{\n-2}}{{\left(\sqrt{t^2+1}\,\right)}^{\n+2-\a}}\,dt
<\int_{-a}^{c_0}
\frac{(t-\frac1a){(t+a)}^{\n-2}}{{\left(\sqrt{t^2+1}\,\right)}^{\n+2-\a}}\,dt
<0,
$$
where the second inequality follows from Lemma \ref{lemma_c_a_m>2} as $a\ge\psi_0(c_0)c_0$. 
It implies $F(\n,\a,a,\xi_0)<0$, which completes the proof. 
\end{proof}
%
%
\begin{lemma}\label{keylemma2}
Suppose $\n\ge2$. 
There is a function $f\colon(1,\n+1)\to\RR_+$ with the following property. 

Suppose $\widetilde\Om$ is a convex set of $\RR^2_{\ge0}=\{(x_1, x_2)\,|\,x_2\ge0\}$ with a non-empty intersection $\widetilde\Ga_0$ with the $x_1$-axis. 
Put $\Om_\d=\big(\widetilde\Om+\d D^2\big)\cap\RR^2_{\ge0}$ and let $\Ga_\d$ be the closure of the intersection of $\partial\Om_\d$ and the upper half plane {\rm (Figure \ref{parallel_body})}. 
Then if $1<\a<\n+1$ and if $\d\ge f(\a)\cdot\textrm{\rm diam}(\widetilde\Om)$ then for any point $x=(\xi,0)\in\widetilde\Ga_0$ we have 
\[
\int_{\Ga_\d}{|x-y|}^{\an-2}(\xi-y_1)\,{y_2}^{\n-2}\,dy_2<0\,. 
\]
\end{lemma}

\begin{proof}
Let us write $\Om=\Om_\d$ and $\Ga=\Ga_\d$ in what follows. 

Suppose $x\in \widetilde\Ga_0$. 
If $\Om''$ is a compact domain that does not contain $x$ then 
\[\begin{array}{l}
\displaystyle \int_{\partial\Om''}{|x-y|}^{\an-2}(\xi-y_1)\,{y_2}^{\n-2}\,dy_2\\[4mm]
\displaystyle =\int_{\Om''}{|x-y|}^{\an-4}\left\{(\n+1-\a)(y_1-\xi)^2-{y_2}^2\right\}{y_2}^{\n-2}\,dy_1dy_2
\end{array}\]
Note that the integrand above is positive if $|y_2|<\sqrt{\n+1-\a\,}\,|y_1-\xi|$ and negative if $|y_2|>\sqrt{\n+1-\a\,}\,|y_1-\xi|$. 

Two lines $y_2=\pm\sqrt{\n+1-\a\,}\,(y_1-\xi)$ intersect $\Ga$ in a point each as $\Om$ is convex. 
Let $L$ be a line through the two intersection points. 
Remark that the line $L$  does not have any other intersection points with $\partial\Om$ as $\Om$ is convex. 
\begin{figure}[htbp]
\begin{center}
\begin{minipage}{.42\linewidth}
\begin{center}
\includegraphics[width=\linewidth]{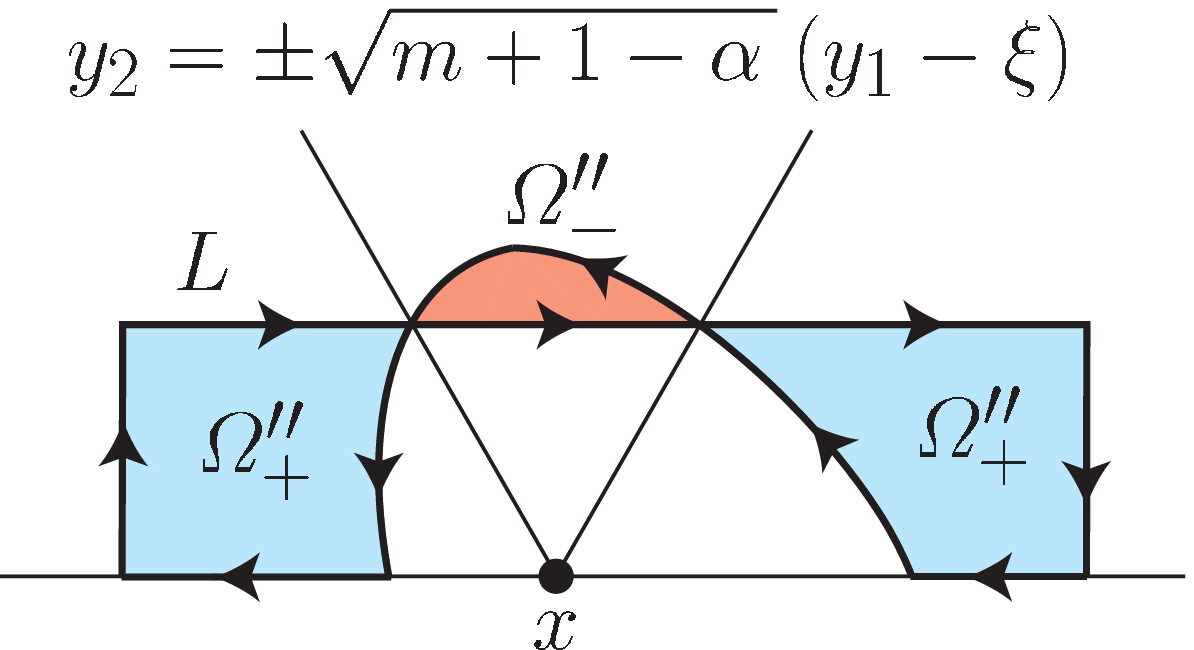}
\caption{Domain $\Om'$ (the rectangle) when $L$ is parallel to the $y_1$-axis. The arrows indicate the orientation of $\Om''=\Om-\Om'$ }
\label{new_domain_3bis}
\end{center}
\end{minipage}
\hskip 0.4cm
\begin{minipage}{.48\linewidth}
\begin{center}
\includegraphics[width=\linewidth]{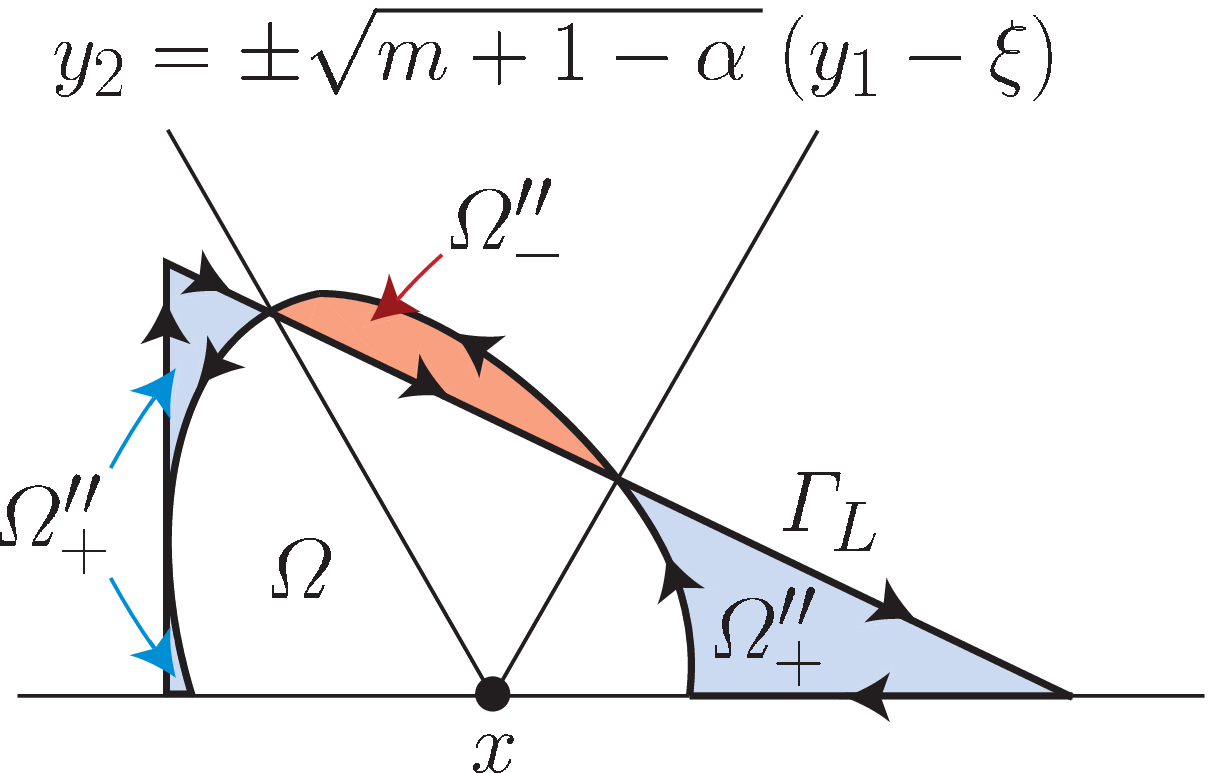}
\caption{Domain $\Om'$ (the triangle) when $L$ is not parallel to the $y_1$-axis. }
\label{new_domain_1bis}
\end{center}
\end{minipage}
\end{center}
\end{figure}
We construct a new domain $\Om'$, which is a rectangle or a triangle, according to whether $L$ is parallel to the $y_1$-axis or not, bounded by a line segment of $L$ (denoted by $\Ga_L$), a line segment of the $y_1$-axis, and some vertical line segments as is indicated in figures \ref{new_domain_1bis} and \ref{new_domain_3bis}. 
When $\Om'$ is a triangle, we take the vertical line segment as close to the point $x$ as possible. 

Put $\Om''_+=\Om'\setminus(\Om'\cap\Om)$ and $\Om''_-=\Om\setminus(\Om'\cap\Om)$ (domains in blue (light gray) and red (dark gray) respectively in figures \ref{new_domain_1bis} and \ref{new_domain_3bis}), and $\Om''=\Om-\Om'$. 
Note that $\Om''=(-\Om''_+)\cup\Om''_-$. 
Then, 
\begin{equation}\label{f_keylemma2-1}
\begin{array}{l}
\displaystyle \int_\Ga{|x-y|}^{\an-2}(\xi-y_1)\,{y_2}^{\n-2}\,dy_2
=\displaystyle \int_{\partial\Om}{|x-y|}^{\an-2}(\xi-y_1){y_2}^{\n-2}\,dy_2\\[4mm]
=\displaystyle \int_{\partial\Om''}{|x-y|}^{\an-2}(\xi-y_1){y_2}^{\n-2}\,dy_2
+\int_{\partial\Om'}{|x-y|}^{\an-2}(\xi-y_1)\,{y_2}^{\n-2}\,dy_2. 
\end{array}
\end{equation}

As 
$$
\stackrel{\!\!\!\circ\,{}}{\Om''_+}\subset \left\{(y_1,y_2)\,|\,|y_2|<\sqrt{\n+1-\a\,}\,|y_1|\right\}, \>\>
\stackrel{\!\!\!\circ\,{}}{\Om''_-}\subset \left\{(y_1,y_2)\,|\,|y_2|>\sqrt{\n+1-\a\,}\,|y_1|\right\},
$$
%
the first term of the right hand side of \eqref{f_keylemma2-1} satisfies 
\[\begin{array}{l}
\displaystyle \int_{\partial\Om''}{|x-y|}^{\an-2}(\xi-y_1)\,{y_2}^{\n-2}\,dy_2\\[4mm]
=\displaystyle -\int_{\Om''_+}{|x-y|}^{\an-4}\left\{(\n+1-\a)(y_1-\xi)^2-{y_2}^2\right\}\,{y_2}^{\n-2}\,dy_1dy_2\\[4mm]
\displaystyle \phantom{=}+\int_{\Om''_-}{|x-y|}^{\an-4}\left\{(\n+1-\a)(y_1-\xi)^2-{y_2}^2\right\}\,{y_2}^{\n-2}\,dy_1dy_2\\[4mm]
<0. 
\end{array}\]

The second term of the right hand side of \eqref{f_keylemma2-1} can be estimated as follows. 
Notice that $\partial\Om'$ consists of a line segment of $L$, which we denote by $\Ga_L$, vertical edges, which we denote by $\Ga_v$, and a horizontal edge on the $y_1$-axis, where the integral vanishes. 
As the orientation of $\Ga_v$ is upward on the right edge and downward on the left edge, we have 
$$\displaystyle \int_{\Ga_v}{|x-y|}^{\an-2}(\xi-y_1)\,{y_2}^{\n-2}\,dy_2<0.$$ 
Therefore, it remains to show 
\[\int_{\Ga_L}{|x-y|}^{\an-2}(\xi-y_1)\,{y_2}^{\n-2}\,dy_2<0 \]
when $\Ga_L$ is not parallel to the $y_1$-axis. It is equivalent to show that $F(\n,\a,a,\xi)<0$ if $a$ and $\xi$ satisfy some conditions which are derived from the condition for $\d$. 

\smallskip
We may assume without loss of generality that the slope of $\Ga_L$ is negative. 
Put $\diam=\textrm{diam}(\widetilde\Om)$. 
Let $z$ (or $w$) be the intersection point of the $y_1$-axis and $\Ga_v$ (or $\Ga_L$ respectively). 
Let $p$ and $q$ be the intersection points of $\Ga_L$ and the lines through $x$ with slopes $\pm\sqrt{\n+1-\a}$ (Figure \ref{new_domain_1bis3}). 
Then 
\begin{figure}[htbp]
\begin{center}
\includegraphics[width=.35\linewidth]{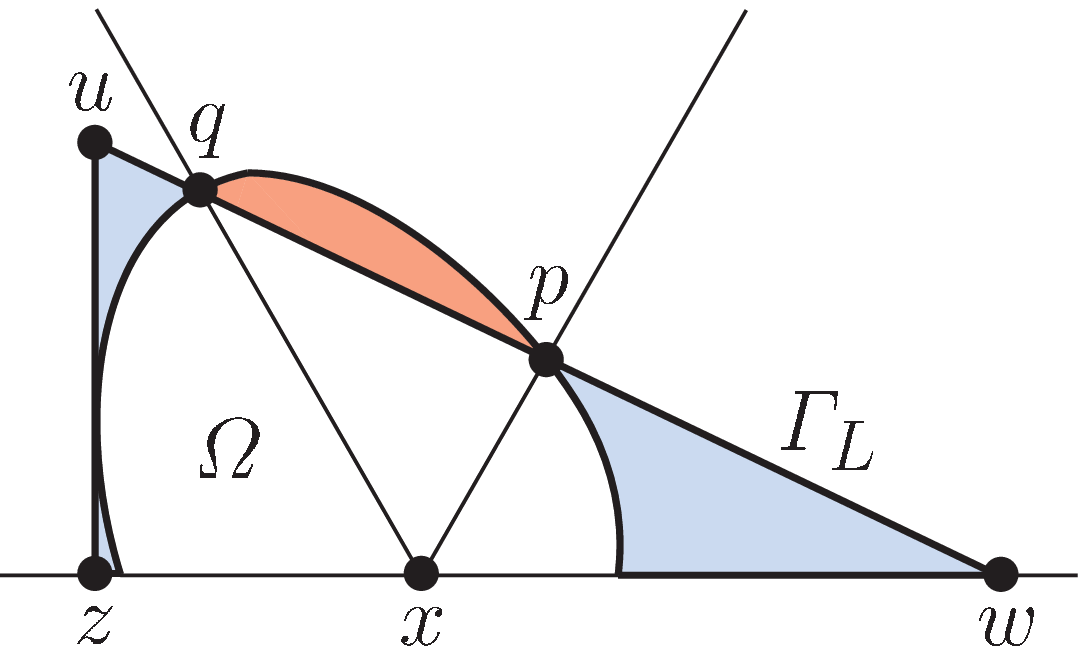}
\caption{}
\label{new_domain_1bis3}
\end{center}
\end{figure}
$$
\d\le|x-p|,|x-q|\le\d + d.
$$
Therefore, the slope of $\Ga_L$ is not greater than $\frac{d}{2\d+d}\sqrt{\n+1-\a}$, and hence 
\begin{equation}\label{x-w}
|x-w| \ge \frac\d{\sqrt{\n+2-\a}}+\frac{\sqrt{\n+1-\a}}{\sqrt{\n+2-\a}}\,\d\,\frac{2\d+d}{d\sqrt{\n+1-\a}}
=\frac{2\d(\d+d)}{d\sqrt{\n+2-\a}}.
\end{equation}
On the other hand, as we take $\Ga_v$ as close to $x$ as poosible, we have 
\begin{equation}\label{x-z}
|x-z|\le\d+d.
\end{equation}
%
%

\smallskip
(i) Suppose $\n=2$. Put $f(\a)=\frac12\sqrt{4-\a}$. 
Then, if $\d\ge f(\a)d$ then 
$$
|x-w|\ge\frac{2\d(\d+d)}{d\sqrt{4-\a}}\ge \d+d \ge|x-z|. 
$$
Lemma \ref{key_lemma_m=2} implies $\int_{\Ga_L}{|x-y|}^{\a-4}(\xi-y_1)\,dy_2<0.$ 

\smallskip
(ii) Suppose $\n\ge3$. Let $u$ be the intersection point of $\Ga_L$ and $\Ga_v$. 
Then $|u-z|\ge\d\sqrt{\frac{\n+1-\a}{\n+2-\a}}$. 
Put 
$$
\xi_0=2\sqrt{\frac{\n+2-\a}{\n+1-\a}}. 
$$
If we assume $\frac\d{d}\ge1$ then \eqref{x-w} and \eqref{x-z} imply
\begin{equation}\label{f_m>2_2}
\frac{|x-z|}{|u-z|}\le\frac{\d+d}{\d}\sqrt{\frac{\n+2-\a}{\n+1-\a}}
\le \xi_0. 
\end{equation}
On the other hand, as the the slope of $\Ga_L$ is not greater than $\frac{d}{2\d+d}\sqrt{\n+1-\a}$ we have 
\begin{equation}\label{w-z/u-z}
\frac{|w-z|}{|u-z|}\ge\frac{2(\d+d)}{d\sqrt{\n+1-\a}}
=2\frac{1+\frac\d{d}}{\sqrt{\n+1-\a}}. 
\end{equation}
Put 
\begin{equation}\label{def_f_m>2}\begin{array}{rcl}
f(\a)&=&\displaystyle \frac{\sqrt{\n+1-\a}}2\>\psi_\a\left(\frac{2\xi_0^2+1}{\xi_0}\right)\frac{2\xi_0^2+1}{\xi_0}-1\\[4mm]
&=&\displaystyle \frac{\sqrt{\n+1-\a}}2\left(2+\frac3{{\left\{1+{{\left[{4\left(4\sqrt{\frac{\n+2-\a}{\n+1-\a}}+\frac12\sqrt{\frac{\n+1-\a}{\n+2-\a}}\,\right)^2+1}\right]}^{-\frac{(\n+2-\a)}2}}\right\}}^{\frac1{\n-2}}-1}\right)\\[14mm]
&&\displaystyle \times\left(4\sqrt{\frac{\n+2-\a}{\n+1-\a}}+\frac12\sqrt{\frac{\n+1-\a}{\n+2-\a}}\,\right)
-1.
\end{array}
\end{equation}
Remark that $f(\a)\ge3$ and hence if $\frac\d{d}\ge f(\a)$ then the assumption $\frac\d{d}\ge1$ above is satisfied. 
%

If $\frac\d{d}\ge f(\a)$ then \eqref{w-z/u-z} implies 
\begin{equation}\label{f_m>2_3}
\frac{|w-z|}{|u-z|}\ge\psi_\a\left(\frac{2\xi_0^2+1}{\xi_0}\right)\frac{2\xi_0^2+1}{\xi_0}. 
\end{equation}
Then by \eqref{f_m>2_2} and \eqref{f_m>2_3}, Corollary \ref{cor_lemma_c_a_m>2} implies 
\[\int_{\Ga_L}{|x-y|}^{\an-2}(\xi-y_1)\,{y_2}^{\n-2}\,dy_2<0. \]
\end{proof}

\begin{theorem}
\label{main_thm}
Suppose $\n\ge2$ and $1<\a<\n+1$. 
Let $\widetilde\Om$ be a compact convex set in $\RR^\n$ with a piecewise $C^1$ boudnary. 
If $\d\ge f(\a)\cdot \textrm{\rm diam}(\widetilde\Om)$, where $f(\a)$ is given in Lemma \ref{keylemma2}, then \eqref{f_second_der} holds on $\widetilde\Om$ for any $j$ $(1\le j\le \n)$. 
%
\end{theorem}

\begin{proof}
Put $\Om=\widetilde\Om+\d B^\n$ in what follows. 
Suppose $x\in\widetilde\Om$. 
By the symmetry, we may assume that $j=1$ and that $x$ is on the $x_1$-axis. 
We omit the proof for the case when $\a=\n$ as it is same as that for the case when $1<\a<\n$. 

(i) The case when $\n=2$. 
Recall \eqref{f_second_partial_derivative_boundary}: 
\[
\frac{\partial^2 V_{\Om}^{(\a)}}{\partial x_1{}^2}(x)
=(2-\a)\int_{\partial\Om}{|x-y|}^{\a-4}(x_1-y_1)\,dy_2.
\]
Divide $\partial\Om$ into two parts by the $x_1$-axis, and lemma \ref{keylemma2} implies the conclusion. 

\smallskip
(ii) The case when $\n\ge3$. 
We use the orthogonal decomposition 
\[\RR^\n=\RR\oplus\RR^{\n-1}=\langle x_1\rangle\oplus\langle x_2,\ldots,x_m\rangle.\]
Suppose the intersection of $\Om$ and the $x_1$-axis is given by $[x_1^{min}, x_1^{max}]$. 

Let $S^{\n-2}$ be the unit sphere in $\RR^{\n-1}$. 
Suppose $\theta_2,\ldots,\theta_{\n-1}$ are local coordinates of $S^{\n-2}$. 
Put $\theta=(\theta_2,\ldots,\theta_{\n-1})$, and let $\gamma(\theta)$ be the corresponding point on $S^{\n-2}$. 
%
Let $\varPi_{\gamma(\theta)}$ be a half $2$-plane in $\RR^\n$ with the axis being the $x_1$-axis that contains the point $\gamma(\theta)$. 

Assume that $\partial\Om$ can locally be parametrized by 
\[
\Phi(t,\theta)=\left(f(t,\theta), \, g(t,\theta)\gamma(\theta)\right)\in\RR\oplus\RR^{\n-1} \hspace{0.5cm} (t_0(\theta)\le t\le t_1(\theta))
\]
so that the following conditions are satisfied. 
\begin{itemize}
\item $f$ and $g$ are piecewise $C^1$-functions with $(f_t)^2+(g_t)^2>0$, 
\item $f(t_0(\theta),\theta)=x_1^{max}, \>f(t_1(\theta),\theta)=x_1^{min}$, 
\item $g(t,\theta)\ge0$, namely $\Phi(t,\theta)\in \varPi_{\gamma(\theta)}$, and $g(t_0(\theta),\theta)=g(t_1(\theta),\theta)=0$, 
%
\end{itemize}
Then, if we put $\Ga_{\gamma(\theta)}=\partial\Om\cap \varPi_{\gamma(\theta)}$ then $\Ga_{\gamma(\theta)}$ can be expressed with respect to the $x_1$-axis and an orthogonal axis in $\varPi_{\gamma(\theta)}$ by (Figure \ref{Gamma_theta})
\[
\bar y(t,\theta)=(\bar y_1, \bar y_2)=
\left(f(t,\theta),g(t,\theta)\right) \hspace{0.5cm} (t_0(\theta)\le t\le t_1(\theta)).
\]
%
%
\begin{figure}[htbp]
\begin{center}
\includegraphics[width=.4\linewidth]{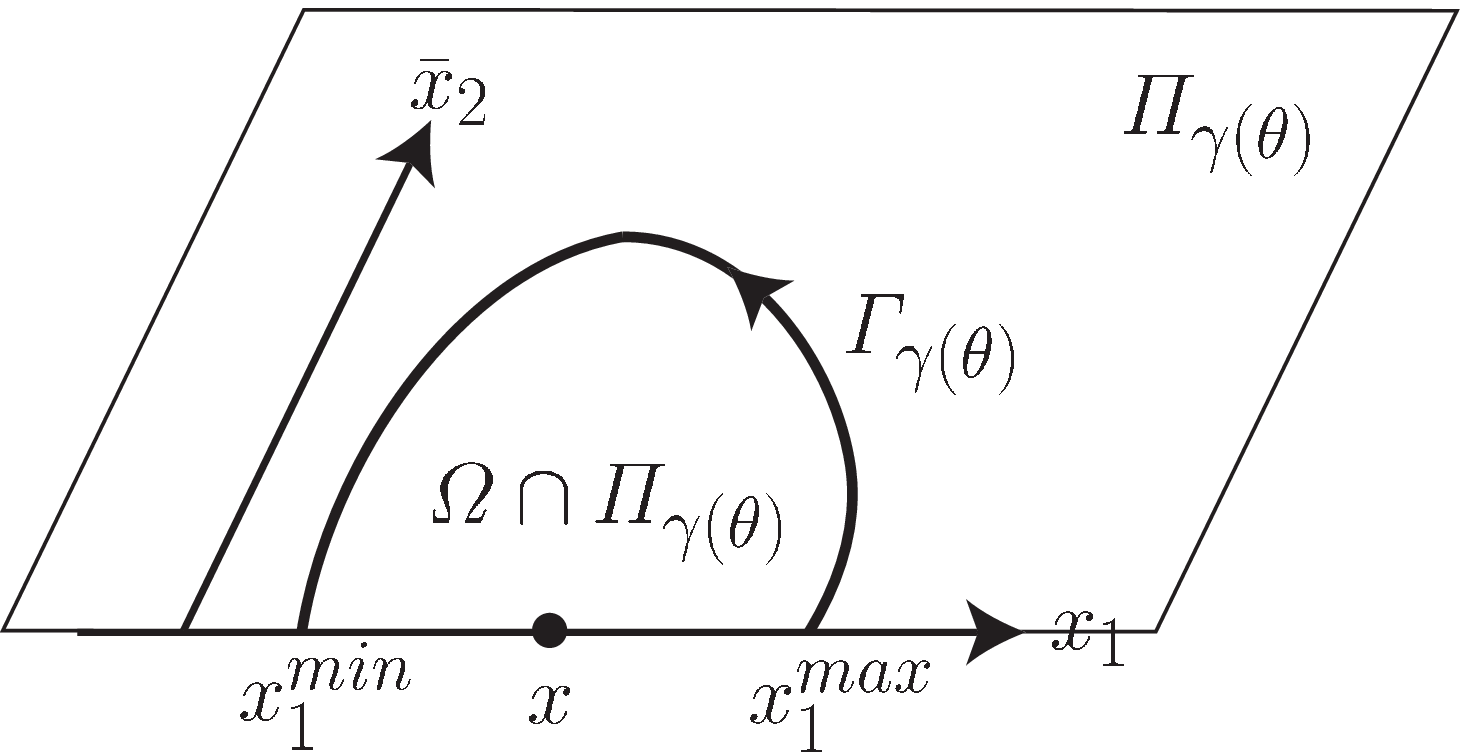}
\caption{}
\label{Gamma_theta}
\end{center}
\end{figure}
Put 
\[\nu=\frac{\partial\Phi}{\partial t}\times \frac{\partial\Phi}{\partial \theta_2}\times \cdots \frac{\partial\Phi}{\partial \theta_{\n-1}},\]
which is a normal vector to $\partial\Om$. 
Then $\nu$ is an outer normal vector if and only if 
\begin{eqnarray}\label{f_outer_normal_nu}
\left(\Phi(t,\theta)-p_0\right)\cdot\nu=\left|\Phi(t,\theta)-p_0 \>\> \frac{\partial\Phi}{\partial t} \>\> \frac{\partial\Phi}{\partial \theta_2} \>  \cdots  \> \frac{\partial\Phi}{\partial \theta_{\n-1}}\right|>0
\end{eqnarray}
for any point $p_0$ in $\interior\Om$ as $\Om$ is convex. 
When $f_t\ne0$ we can take $p_0$ in $\varPi_{\gamma(\theta)}$ so that $p_0$ has the same $x_1$-coordinate as $\Phi(t,\theta)$. 
Then $\Phi(t,\theta)-p_0$ is a positive multiple of $\left(0, -(\textrm{sgn}\, f_t)\gamma(\theta)\right)$. 
Therefore, if $f_t\ne0$ then \eqref{f_outer_normal_nu} is equivalent to 
\[\begin{array}{rcl}
0&<&\displaystyle 
\left|\begin{array}{ccccc}
0 & f_t & g_{\theta_2} & \cdots & g_{\theta_{\n-1}}\\
-(\textrm{sgn}\, f_t)\gamma \>&\> g_t\gamma \>&\> g_{\theta_2}\gamma+g\gamma_{\theta_2} \>& \cdots &\> g_{\theta_{\n-1}}\gamma+g\gamma_{\theta_{\n-1}}
\end{array}\right| \\[4mm]
&=&\displaystyle g^{\n-2}|f_t|
\left|\begin{array}{ccccc}
1 & 0 & 0 & \cdots & 0 \\
0 \>&\> \gamma \>&\> \gamma_{\theta_2} \>& \cdots &\> \gamma_{\theta_{\n-1}}
\end{array}\right| \\[4mm]
&=&\displaystyle g^{\n-2}|f_t|\,\left|\gamma \>\> \gamma_{\theta_2} \> \cdots \> \gamma_{\theta_{\n-1}}\right|. 
\end{array}\]
Assume $\theta_2,\ldots,\theta_{\n-1}$ are positive local coordinates of $S^{\n-2}$, i.e. $\left|\gamma \>\> \gamma_{\theta_2} \> \cdots \> \gamma_{\theta_{\n-1}}\right|>0$. 
Then $\nu$ is an outer normal vector to $\partial\Om$. 
This holds even when $f_t=0$ because in this case $\nu$ is outer normal if and only if $(\textrm{sgn}\,g_t)e_1\cdot\nu>0$, which follows from \eqref{e_1_dot_nu} below. 

\smallskip
On the other hand, 
\begin{eqnarray}\label{e_1_dot_nu}
\begin{array}{rcl}
e_1\cdot\nu&=&\displaystyle 
\left|\begin{array}{ccccc}
1 & f_t & g_{\theta_2} & \cdots & g_{\theta_{\n-1}}\\
0 \>&\> g_t\gamma \>&\> g_{\theta_2}\gamma+g\gamma_{\theta_2} \>& \cdots &\> g_{\theta_{\n-1}}\gamma+g\gamma_{\theta_{\n-1}}
\end{array}\right| \\[4mm]
&=&\displaystyle g^{\n-2}g_t\left|\gamma \>\> \gamma_{\theta_2} \> \cdots \> \gamma_{\theta_{\n-1}}\right|. 
\end{array}
\end{eqnarray}
Since 
\[
\begin{array}{rcl}
\displaystyle dS^{\n-2}&=&\displaystyle \left|\gamma \>\> \gamma_{\theta_2}\>\> \cdots \>\> \gamma_{\theta_{\n-1}}\right|\,d\theta_2\cdots d\theta_{\n-1},\\[2mm]
d\sigma
&=&|\nu|\,dt\,d\theta_2\cdots d\theta_{\n-1},\\[2mm]
n&=&\nu/|\nu|,
\end{array}
\]
we have 
\[
\begin{array}{rcl}
e_1\cdot n\,\,d\sigma&=&\displaystyle g^{\n-2}g_t\left|\gamma \>\> \gamma_{\theta_2} \> \cdots \> \gamma_{\theta_{\n-1}}\right|\,dt\,d\theta_2\cdots d\theta_{\n-1}\\[2mm]
&=&\displaystyle g^{\n-2}g_t\,dt\,dS^{\n-2}.
\end{array}
\]

Therefore, \eqref{f_second_partial_derivative_boundary} implies that 
\[\begin{array}{rcl}
\displaystyle \frac{\partial^2 V_{\Om}^{(\a)}}{\partial x_1{}^2}(x)
&=&\displaystyle (\n-\a)\int_{\partial\Om}{|x-y|}^{\an-2}(x_1-y_1)\,e_1\cdot n\,d\sigma (y)\\[4mm]
&=&\displaystyle (\n-\a)\int_{S^{\n-2}}\left(\int_{t_0(\theta)}^{t_1(\theta)}|x-\bar y|^{\an-2}(x_1-\bar y_1)g^{\n-2}g_{t}\,dt\right)dS^{\n-2}.
\end{array}\]
By lemma \ref{keylemma2} 
$$
\int_{t_0(\theta)}^{t_1(\theta)}|x-\bar y|^{\an-2}(x_1-\bar y_1)g^{\n-2}g_{t}\,dt
=\int_{\Ga_{\gamma(\theta)}}|x-\bar y|^{\an-2}(x_1-\bar y_1)\,{\bar y_2}{}^{\!\n-2}\,d\bar y_2
<0
$$
%
for each point $\gamma(\theta)$ in $S^{\n-2}$, 
which completes the proof. 
\end{proof}
\begin{corollary}
Suppose $\n\ge2$ and $1<\a<\n+1$. 
For any compact convex set $\widetilde\Om$ in $\RR^\n$ with a piecewise $C^1$ boundary, if $\d\ge f(\a)\cdot \textrm{\rm diam}(\widetilde\Om)$, where $f(\a)$ is given in Lemma \ref{keylemma2}, then $\widetilde\Om+\d B^\n$ has a unique $r^\an$-center. 
\end{corollary}

When $\n=2$ we have $\sup_{1<\a<3}f(\a)=\sqrt3$, so if we put $\varphi(2)=\sqrt3$ we completes the proof of Theorem \ref{main_thm1} for the case when $\n=2$. 

\medskip
When $\n\ge3$, unfortunately we have $\sup_{1<\a<\n+1}f(\a)=+\infty$ as $\displaystyle \lim_{\a\nearrow\n+1}f(\a)=+\infty$. 

\begin{lemma}\label{last_lemma}
Suppose $\n\ge3$. 
For any $b>0$ there is $\a_0=\a_0(b)$ with $\n<\a_0<\n+1$ such that for any compact convex set $\widetilde\Om$ in $\RR^\n$ with a piecewise $C^1$ boundary, if $\d\ge b\cdot \textrm{\rm diam}(\widetilde\Om)$ then $\widetilde\Om+\d B^\n$ has a unique $r^\an$-center if $\a_0\le\a<\n+1$. 
\end{lemma}
\begin{proof}
Suppose $\widetilde\Om$ has diameter $\diam$ and $x\in\widetilde\Om$. 
Let $C_j(\a)$ be the cone with vertex $x$ given by 
$$
C_j(\a)=\Big\{y\,\Big|\,-(\n+1-\a)(x_j-y_j)^2+\sum_{i\ne j}(x_i-y_i)^2\le0\Big\}. 
$$
The radial function of $\widetilde\Om+\d B^\n$ with respect to $x$ defined by $\rho(v)=\sup\{t\ge0\,|\,x+tv\in\widetilde\Om+\d B^\n\}$ $(v\in S^{\n-1})$ satisfies $\delta\le\rho(v)\le\delta+d$ for any $v$. 
Therefore 
\begin{equation}\label{final_f_array}
\begin{array}{l}
\displaystyle \frac1{\a-\n}\cdot\frac{\partial^2 V_{\widetilde\Om+\d B^\n}^{(\a)}}{\partial x_j{}^2}(x)
=\displaystyle \int_{\widetilde\Om+\d B^\n}{|x-y|}^{\an-4}\left(-(\n+1-\a)(x_j-y_j)^2+\sum_{i\ne j}(x_i-y_i)^2\right)\!d\mu (y)\\[6mm]
\ge\displaystyle \int_{B_{\d+d}^\n(x)\cap C_j(\a)}{|x-y|}^{\an-4}\left(-(\n+1-\a)(x_j-y_j)^2+\sum_{i\ne j}(x_i-y_i)^2\right)\!d\mu (y)\\[6mm]
\displaystyle +\int_{B_{\d}^\n(x)\cap {C_j(\a)}^c}{|x-y|}^{\an-4}\left(-(\n+1-\a)(x_j-y_j)^2+\sum_{i\ne j}(x_i-y_i)^2\right)\!d\mu (y).
\end{array}
\end{equation}

Define $g\colon(\n,\n+1)\times\RR_+\to\RR$ by 
$$
g(\a,\beta)=\int_{X_{\a,\beta}}\frac{-(\n+1-\a){y_1}^2+\sum_{i>1}(x_i-y_i)^2}{{|y|}^{\n+4-\a}}\,d\mu (y),
$$
where 
$$
X_{\a,\beta}=B^\n\cup\left(B_{1+\frac1\beta}^\n\cap C_1(\a)\right).
$$
Remark that $g(\a,\beta)$ is an increasing function of $\beta$. 
Fix $b>0$. 
As $g(\a,b)$ is continuous with respect to $\a$ and $g(\n+1,b)>0$, there is $\a_0\in(\n,\n+1)$ such that if $\a_0\le\a<\n+1$ then $g(\a,b)>0$, which completes the proof as the right hand side of \eqref{final_f_array} is proportional to $g(\a,b)$. 
\end{proof}

Suppose $\n\ge3$. 
Put 
$$
\varphi(\n)=\max\left\{10, \sup_{1<\a\le\a_0(10)}f(\a)\right\}
=\max\left\{10, \max_{1\le\a\le\a_0(10)}f(\a)\right\},
$$
where $f(\a)$ is given by \eqref{def_f_m>2} and $\a_0$ is given in Lemma \ref{last_lemma}. 
Then, Theorem \ref{main_thm} and Lemma \ref{last_lemma} implies Theorem \ref{main_thm1} for the case when $\n\ge3$. 

\begin{remark}\rm 
In \cite{O}, using the same renormalization process of defining energy functionals of knots (\cite{O1}, \cite{O2}), we renormalized $V_\Om^{(\a)}(x)$ so that it is well-defined for $\a\le0$ and $x\in\interior\Om$. 
The $r^\an$-center of $\Om$ for $\a\le0$ can be defined in a similar way: it is a point in $\interior\Om$ where $V_\Om^{(\a)}\big|_{\interior\Om}$ attains the maximum value. 
We can show $\displaystyle \frac{\partial^2 V_{\Om}^{(\a)}}{\partial {x_j}^2}<0$ on $\interior\Om$ for any $j$ if $\a\le1$ and $\Om$ is convex by a similar way as in Lemma \ref{keylemma2} and Theorem \ref{main_thm}. This gives an alternative proof of the uniqueness of the $r^\an$-center of $\Om$ when $\a\le1$ and $\Om$ is convex. 
\end{remark}

Department of Mathematics and Information Sciences, 

Tokyo Metropolitan University, 

1-1 Minami-Ohsawa, Hachiouji-Shi, Tokyo 192-0397, JAPAN. 

E-mail: ohara@tmu.ac.jp


\begin{thebibliography}{E} 






\bibitem[GNN]{GNN}B.~Gidas, W. M. ~Ni, and L.~Nirenberg, {\em Symmetry and related properties via the maximum principle,} Comm. Math. Phys. {\bf 68} (1979), 209\,--\,243.

\bibitem[H]{H}I.~Herburt, {\em On the Uniqueness of Gravitational Centre,}  Mathematical Physics, Analysis and Geometry {\bf 10} (2007), 251\,--\,259



\bibitem[L1]{Lu75}E.~Lutwak, {\em Dual mixed volumes,} Pacific J. Math. {\bf 58} (1975), 531\,--\,538. 

\bibitem[L2]{Lu88}E.~Lutwak, {\em Intersection bodies and dual mixed volumes,} Advances in Mathematics {\bf 71} (1988), 232\,--\,261.



\bibitem[M]{M1}M.~Moszy\'nska, {\em Looking for selectors of star bodies,} Geom. Dedicata {\bf 81} (2000), 131\,--\,147.


\bibitem[O1]{O1}J.~O'Hara, {\em Energy of a knot,} Topology {\bf 30} (1991), 241\,--\,247.

\bibitem[O2]{O2}J.~O'Hara, {\em  Family of energy functionals of knots,} 
Topology Appl. {\bf 48} (1992), 147\,--\,161.



\bibitem[O3]{O}J.~O'Hara, {\em Renormalization of potentials and generalized centers}, to appear in Adv. Appl. Math., available at {arXiv:1008.2731}. 







\end{thebibliography}
\end{document}